\documentclass{article}

\usepackage{amsmath,amsfonts,amssymb}
\usepackage[symbol]{footmisc}

\begin{document}

\title{Refining Calculus Pedagogy}

\author{Parthasarathy Srinivasan}
\date{}

\maketitle

\begin{abstract}
\noindent There have been several modifications of how basic
calculus has been taught, but very few of these modifications have
considered the computational tools available at our disposal.
Here, we present a few tools that are easy to develop and use.
Doing so also addresses a different way to view calculus, and
attempts to fill the gaps in students' understanding of both
differentiation and integration. We will describe the basics of
both these topics in a way that might be much more useful and
relevant to students, and hence possible ways in refining calculus
pedagogy to make calculus more accessible to them. For
integration, an elementary development of Gaussian quadrature
using basic linear algebra is presented. This numerical method can
be extended to integrate functions over various domains in higher
dimensions, a subject that is not currently well covered in
multivariable calculus or other Mathematics courses. We also
briefly discuss series that may be more useful than some of those
taught in current calculus courses.
\end{abstract}

\subsection*{Introduction}

\noindent Basic calculus has two main ideas -- Differentiation and
integration. Currently, the introduction of differentiation to
students is done via limits, and integration is mainly introduced
by using certain techniques of computing antiderivatives. However,
the geometric aspect of the derivative is not retained by most
students, and the issue of computing integrals when
antiderivatives cannot be obtained using standard methods is
ignored. In this paper, we address both of these issues, and
emphasize the role that computers can play in enhancing how
students can get a better understand integration. A common theme
is a scientific process of learning, where students can verify
their work, which may help them eliminate the mistakes they are
making in their computations.

\subsection*{Differentiation}

\noindent Finding the derivative of a function $y = f(x)$ of a
single variable at a point on its graph is almost trivial, as it
is just a local estimate of the slope of the tangent line at this
point. Suppose this point on the graph of the function is $(a,
f(a))$ where the derivative exists. The slope of the tangent line
can be estimated efficiently by a wise choice of points, typically
$(a-h, f(a-h))$ and $(a+h, f(a+h))$ for small values of $h$, and
then computing the slope of the secant line between these two
points. This secant estimate for the tangent is very reasonable
for most practical applications. The slope of this secant line can
also be used to approximate the tangent line for the function at
$(a, f(a))$, thus locally approximating the function by a line at
that point. In fact, the slope using the secant approximation
$(f(x+h) - f(x))/h$ is built into the definition of the
derivative. Once differentiability of a function has been
established, the above description requires the computation of at
most two function values $f(a+h)$ and $f(a-h)$ for as small a
value of $h$ as computationally feasible, and this depends on the
accuracy of the estimate required. Can we do better than this
computationally? Yes, we can by computing $f'(a)$, as it requires
the evaluation of just a single function value. However, this
requires the additional computation of $f'(x)$, which could be
somewhat involved depending on the nature of $f(x)$. While it may
not be pedagogically relevant, the issue of whether $f'(a)$ is
exact is debatable. For instance, if $f'(a) = 2e^2$ for some
function $f$ and some point $a$, then we can only evaluate this to
a desired number of digits of accuracy, which is not very
different from the numerical method.  It can be pointed out here
to students that by computing $(f(x+h) - f(x))/h$ for positive and
negative values of $h$ that functions with jumps don't have a
derivative at the jump, and functions like $f(x) = |x|$ don't have
a derivative at $x=0$.
\\

Pedagogically, once the functional form of a derivative  $f'(x)$
is computed, the numerical method described above can be used by
students to verify if they have the correct functional form of the
derivative by computing $f'(a)$ and comparing it to the secant
approximation at that point for a particular value of $a$ and a
small value of $h$. For instance, if $f(x) = (x - 2)/(x^2 + 4), a
= 2$, and $h = 10^{-4}$, then the difference between $f'(a)$ and
the secant approximation is around $1.6 \times 10^{-10}$. This
gives students validation of their analytical work. If their
answers don't match, then it invites the student to go back and
determine where they might have made errors. If there are any
doubts, then students can always compute this at a different value
of $a$ to further check their work.  This should be a routine
aspect of the course when students are taught derivatives using
the various rules. \\

It must be noted here that starting with the notion of a limit and
making it somewhat abstract makes the subject much more
challenging to students than it actually is, but this may
demystify it. It may also be useful to introduce limits
using this idea. \\

In higher dimensions for a real-valued function $y =
f(x_1,x_2,\ldots,x_n)$, one generalization of the derivative is
the directional derivative in the direction of a particular
vector. Again, once differentiability has been established, the
directional derivatives can be computed using the partial
derivatives using their secant approximations, and this requires
the computation of at least $n+1$ function values $f(a_1 + h_1,
a_2, \ldots, a_n), f(a_1, a_2 + h_2, \ldots, a_n), \ldots, f(a_1,
a_2, \ldots, a_n + h_n)$ for values of $h_1, h_2, \ldots, h_n$
(typically all equal), as well as $f(a_1, a_2, \ldots, a_n)$. The
analytic form requires evaluating $n$ partial derivatives at the
desired value, but has the additional cost of computing the
various partial derivatives. Other notions of the derivative like
the curl and divergence have different types of complexities
depending on the nature of the function, but all require the
computation of various partial derivatives. \\

As applications for a single variable, the secant method or the
Newton method is used for finding solutions of $f(x) = a$ for some
fixed value of $a$ computationally, but is seldom even introduced
in a calculus course, even though the Newton method is presented
in most calculus textbooks. Using the secant approximation often
for the derivatives might make the use of these computational
methods easier for students. \\

Differentiation can sometimes be used to find the global maximum
or minimum of a smooth function, but the results are not
guaranteed. This is a difficult topic in higher dimensions, but
has many practical applications. It is seldom covered in any
detail in any mathematics course except for a few toy problems in
calculus. Some numerical methods may work in some cases, but they
are in general a little beyond the scope of students in calculus
courses. In any case, there are very few analytic results in this
area. Students must be made aware of this even at the basic
calculus level so that they are aware of the limitations of the
methods presented to them.

\subsection*{Integration}

Integration is non-trivial, as it requires the calculation of the
signed area of a function $f(x)$ bounded by its graph, the
$x$-axis, and the bounds $x=a$ and $x=b$. The question may be
framed as follows - How many computations need to be performed for
say bounded functions with at most finitely many discontinuities
in $[a, b]$ so that we have a good approximation of the integral?
The Fundamental Theorem of Calculus states that whenever we can
obtain an antiderivative $F(x)$ of $f(x)$, then the integral is
$F(b) - F(a)$ ((unlike the case of the derivative, there is no
known general method where the computation at a single point will
suffice other than providing a crude numerical approximation).
This is obviously a game-changer for the computation of the
integral for those functions whose antiderivatives can be
determined using elementary functions. It was clear a long while
ago that this was not possible for all functions, as in the case
of elliptic integrals first studied by Fagnano and Euler. As in
the case of the derivative, it can be debated as to whether this
value is exact or not. \\

If an antiderivative is not easily available or is very difficult
to compute, then we don't address the issue of how to proceed
properly in evaluating an integral in a general calculus course,
even though this should be the primary goal of the integration
part of elementary calculus. Currently, students are only taught
how to compute a few antiderivatives that may or may not be useful
in practice. But the question is better framed by asking how many
computations need to be done in order to get a good approximation
of the integral when the antiderivative cannot be easily computed.
The integrals, for at least the class of functions described above
can be computed numerically in a very efficient way for any
continuous function (and more) using an efficient method like
Gaussian quadrature. We give an elementary description of this
method here. \\

In the discussion below, which is one of the main results
developed in this article, we assume that the antiderivatives for
$x^n$ are known for $n = 0, 1, 2, 3, \ldots$. This can be easily
deduced from formally computing the functional forms of the
derivatives. \\

Now in order approximate $\displaystyle \int_{-1}^1 f(x)\ dx$
using just three points, we have

\[
\int_{-1}^1 f(x)\ dx \approx w_1f(x_1) + w_2f(x_2) + w_3f(x_3).
\]

This can be done in several ways. For instance, the trapezoidal
rule uses points $\{ -1, 0, 1\}$ with corresponding weights
$\{2/3, 4/3, 2/3\}$. \\

We now impose that the integral of any quadratic polynomial using
this method must be exact. Therefore, we have

\begin{align*}
  w_1 + w_2 + w_3 &= \int_{-1}^1 1\ dx = 2 \\
  w_1x_1 + w_2x_2 + w_3x_3 &= \int_{-1}^1 x\ dx = 0 \\
  w_1x_1^2 + w_2x_2^2 + w_3x_3^2 &= \int_{-1}^1 x^2\ dx = 2/3.
\end{align*}

In matrix form, we solve

\begin{equation} \label{eq1}
\left[%
\begin{array}{ccc}
  1 & 1 & 1 \\
  x_1 & x_2 & x_3  \\
  x_1^2 & x_2^2 & x_3^2
\end{array}%
\right] \left[%
\begin{array}{c}
  w_1 \\
  w_2 \\
  w_3
\end{array}%
\right] = \left[%
\begin{array}{c}
  2 \\
  0 \\
  2/3
\end{array}%
\right].
\end{equation}

\medskip

For Gaussian quadrature, we pick $x_j$ so that $P_3(x_j) = 0$ for
$j = 1, 2, 3$, where $P_n(x)$ is the $n^{\textrm{th}}$ degree
Legendre polynomial, and solve for $w_j, j = 1, 2, 3$ in the
system above. \\

The choice of the points $x_j$ is so that $\displaystyle
\int_{-1}^1 x^i P_3(x)\ dx = 0$ for $i = 0, 1, 2$ for our
numerical method. This is because $P_3(x)$ is orthogonal to each
polynomial of degree 2 or lower. Namely if $q(x)$ is any
polynomial of degree 2 or lower, then $\displaystyle \int_{-1}^1
q(x) P_3(x)\ dx = 0$. Now note that from our numerical scheme from
(\ref{eq1}), we have

\[
\int_{-1}^1 x^iP_3(x)\ dx = w_1x_1^iP_3(x_1) + w_2x_2^iP_3(x_2) +
w_3x_3^iP_3(x_3),
\]

\noindent for $i = 0, 1, 2$. This is satisfied by choosing
$P_3(x_j) = 0$ for $j = 1, 2, 3$. Namely, the $x_j's$ are the
zeros of $P_3(x)$, all of which lie in $(-1,1)$. Moreover, let
$p(x)$ be any polynomial of degree less than or equal to five.
Then using the division algorithm, we can write $p(x) = P_3(x)q(x)
+ r(x)$ where both $q(x)$ and $r(x)$ are polynomials that have
degree less than or equal to two. Therefore, $\displaystyle
\int\limits_{-1}^1 p(x)\ dx = \int\limits_{-1}^1 r(x)\ dx$, as
$P_3(x)$ is orthogonal to $q(x)$ and hence $\displaystyle
\int\limits_{-1}^1 P_3(x)q(x)\ dx = 0$. Note that this method

\begin{enumerate}
    \item can be easily extended to $n$ points.
    \item gives the exact integral up to any $2n-1$ degree
    polynomial using $n$ points using a similar argument above. This
    uses the well-known result \cite{Abramowitz1964} that all the roots of
    $P_n(x)$ are distinct and lie in $(-1,1)$, a fact that was known
    to Jacobi \cite{Jacobi1826}.
    \item approximates the integral of any function that is
    well-approximated by any $2n-1$ degree polynomial using $n$
    points. So convergence is rapid.
    \item doesn't need to be restricted to the domain to $[-1,1]$. We can
    use a linear map $u(x) = (b-a)x/2 + (b+a)/2$ to transform $\displaystyle \int_a^b
    f(u)\ du$ for any function $f(u)$ from the function class described above
    on $[a,b]$ to $\displaystyle (b-a)/2 \int_{-1}^1 f\left((b-a)x/2 + (b+a)/2\right)\
    dx$. Note that $(b-a)/2$ is the Jacobian of the
    transformation.
\end{enumerate}

The $n^{\textrm{th}}$ degree Legendre polynomial can be obtained
starting from the vector space spanned by $\{1, x, x^2, \ldots,
x^n \}$ in the interval $[-1,1]$, using the Gram-Schmidt
orthogonalization process to find orthogonal polynomials, and
imposing the condition $P_n(1) = 1$ for all natural numbers $n$.
It is well-known \cite{Abramowitz1964} that all the roots of
$P_n(x)$ are distinct and lie in $(-1,1)$, a fact
that was known to Jacobi \cite{Jacobi1826}. \\

The weights need not be computed by solving systems of equations
like (\ref{eq1}), as there are simple formulas for them
\cite{Hildebrand1956}.

\[
w_i = \frac{2}{(1-x_i^2)[P_n'(x_i)]^2}, \qquad i = 1, 2, \ldots,
n.
\]

Thus, the points and weights in which the function values must be
evaluated need not be computed every time, but can be stored in a
file and used as required. This is similar to tables of logarithms
or exponential functions. In fact, in order to make the
computations more efficient, the values of standard functions may
also be stored and recalled as required. \\

As Gaussian quadrature doesn't use the end points 1 and $-1$ of
the interval of integration, it can be adapted for typical
improper integrals from a calculus course, as long as they
converge. Numerical convergence may be slow due to the unbounded
nature of a function. For instance, $\displaystyle \int_0^1
1/\sqrt{x}\ dx$ with forty quadrature points gives a values of
around 1.9785
(exact value is 2). \\

As with all numerical methods, this can be computationally
expensive depending on the precision of the answer required, but
is very efficient if the number of digits of precision is not very
big (i.e. most practical applications). Note that convergence may
be slow if the function is not well-approximated by polynomials,
like $f(x) = e^{-x^2}$, and a more appropriate choice of basis
might provide faster convergence. \\

Pedagogically, it is important to note that this method can be
easily generalized for integration in higher dimensions of a large
class of functions over several commonly used domains, unlike the
very limited number of integration problems we address in
multivariable calculus courses currently. Results like Green's
Theorem, Stokes' Theorem or the Divergence Theorem can be used to
reduce the dimension of the problem if possible, and this may
reduce the computational complexity of the problem.

\subsection*{Series}

The way series are currently taught could be improved as well. The
various Taylor series developed in calculus courses are usually
not very good in terms of computational efficiency. One potential
use of series is to compute the values of standard functions. This
can be achieved very effectively using shift-and-add algorithms
like CORDIC (coordinate rotation digital computer), first
described by Jack E. Volder \cite{Volder1959}, that avoids the use
of multiplication altogether. The development of this for
computing the values of sine and cosine also uses rotation
matrices from basic linear algebra, and may be a challenge for
students. However, a solid preliminary foundation may make this
easier to understand for them. Though the details are not
presented here, the reader is encouraged to explore this wonderful
topic, and the more general shift-and-add algorithms that
efficiently compute values of various functions.

\subsection*{Discussion}

The approaches described here may be more useful than what is
usually taught in general calculus courses. As mentioned earlier,
it does require some elementary development of linear algebra in
order for students to understand the details of the methods
involved for Gaussian quadrature. This is more an opportunity to
highlight the uses of linear algebra rather than a hinderance. The
Gram-Schmidt process of creating orthogonal vectors starting from
a set of independent vectors and the associated $QR$ decomposition
is a particularly useful and visually appealing ideas in linear
algebra. Various fields where calculus is applied can adapt these
methods to develop their own special ways to use this method more
efficiently. For example, if they suspect that the vector space
spanned by $\{1, x, x^2, \ldots, x^n \}$ may not be good enough
for fast convergence, then they can start with a different basis
suitable to their particular problem to see if convergence is more
rapid. \\

Does this mean that students don't need to know any
antiderivatives? For pedagogy, this is not necessarily true. If
they are asked to just compute an antiderivative $F(x)$ for
$f(x)$, then they have no reasonable way of knowing if their
answer is correct or not. One way to verify is to compute the
derivative of $F(x)$ and check that it is indeed $f(x)$. An
alternate way is to compute $\displaystyle \int_a^b f(x)\ dx =
F(b) - F(a)$ for some choice of $a, b$ where the integral is
well-defined, and compare it to the numerical method described to
determine if they have the correct antiderivative. They can
currently do this with the standard TI-83 or TI-84 calculator as
well, but the difference is that they have some understanding of
how a calculator might
compute this numerical integral. \\

The main reason for introducing students to Gaussian quadrature
over other common methods is that the convergence in this case is
rapid compared to methods like the Midpoint or Trapezoidal rule.
This makes a huge difference when using numerical methods to
compute integrals in two and three dimensions. \\

A note about computers, calculators and their usage. With the
ubiquitous availability of various software including ones
designed for symbolic computations, one can of course use it to
get functional forms of derivatives and many antiderivatives.
However, students will just treat it as a black box without
understanding the details involved. Combining the numerical
methods detailed here along with deriving the functional forms
might go a long way in reducing the difficulty that many students
have with calculus. It will also give them a better understanding
of what the subject is actually about, rather than just some
recipe for pushing symbols around, and which means little for many
of them. Verifying their analytical work along with the numerical
methods is also an important tool for any quantitative subject.
This is seldom done in most courses, and I think calculus
courses are excellent for them to hone this useful skill. \\

We haven't discussed the calculus reform movement that started in
the late 1980's or its subsequent effects here. While there is
nothing new in the material presented here, it is not only
significantly different from that reform movement, but also
different from the basic approach taken in most mathematics
courses. There's talk of conceptual understanding in calculus, but
what exactly does this entail? I fear currently, most students
don't know what either a derivative, an integral, or a series
actually is, and whether they are even useful for something. The
hope with the approach presented here is that students develop
analytical skills along with numerical ones, and are able to apply
and adapt these skills in various fields. \\

Do the details of calculus actually matter? They do, but are only
appealing to someone pursuing a degree in the subject, and much
better covered in a real analysis course, as is done now. There's
enough material available on the subject for the curious student
in any field to explore and learn the details involved. Even in
regular calculus courses, while we spend quite a bit of time going
over the details of a limit, and some details of differentiation,
we largely gloss over these details in integration. As mentioned
earlier, we focus only on computing a few antiderivatives. The
reason for this is that the details for integration are more
involved than for those of limits and derivatives. Ultimately, we
have to decide how much we really want the students in a basic
calculus to know for the ideas introduced in the course to be
useful to them. \\

In summary, smart ideas can save months of unnecessary pedagogy
that's ultimately not very useful. Similar reforms can be done in
other topics like linear algebra and differential equations. This
is just a first step.

\subsection*{Acknowledgements}
I would like to take this opportunity to thank the generous
support of my institution, Cleveland State University. All the
ideas presented here are my own, and do not represent those of
Cleveland State. I would also like to thank Professor Kit Chan
from Bowling Green State University for carefully reading the
manuscript and for helpful suggestions.\\


\subsection*{Disclosure}
No potential competing interest was reported by the author.

\end{document}